# Validity of Borodin & Kostochka Conjecture for a Class of Graphs

Medha Dhurandhar

**Abstract:** **Borodin & Kostochka** conjectured that if $\Delta(G) \geq 9$, then $\chi(G) \leq \max\{\omega, \Delta-1\}$. Here we prove that this Conjecture is true for $\{P_3 \cup K_1\}$-free graphs and $\{K_2 \cup \overline{K_2}\}$-free graphs.

**Introduction:**

In 1977, **Borodin & Kostochka** conjectured that if $\Delta(G) \geq 9$, then $\chi(G) \leq \max\{\omega, \Delta-1\}$ [1]. In 1999, Reed proved the conjecture for $\Delta \geq 10^{14}$ [2]. Also D. W. Cranston and L. Rabern [3] proved it for claw-free graphs and Medha Dhurandhar [4] proved it for $2K_2$-free graphs. Here we prove that the conjecture is true for $\{P_3 \cup K_1\}$-free graphs and $\{K_2 \cup \overline{K_2}\}$-free graphs.

**Notation:** For a graph G, V(G), E(G), $\Delta$, $\omega$, $\chi$ denote the vertex set, edge set, maximum degree, size of a maximum clique, chromatic number of G resply. For $u \in V(G)$, $N(u) = \{v \in V(G) / uv \in E(G)\}$, and $\overline{N(u)} = N(u) \cup (u)$. If $S \subseteq V$, then $<S>$ denotes the subgraph of G induced by S. If C is some coloring of V(G) and if $u \in V(G)$ is colored m in C, then u is called a m-vertex. Also if P is a path in G s.t. vertices on P are alternately colored say i and j, then P is called an i-j path. All graphs considered henceforth are simple and undirected. For terms not defined herein, we refer to Bondy and Murty [5].

Let $H = P_3 \cup K_1$ and $R = K_2 \cup \overline{K_2}$. Note that the only odd cycle in a H-free or R-free graph is $C_5$.

**Main Result 1:** If $\Delta \geq 9$, and G is H-free, then $\chi \leq \max\{\omega, \Delta-1\}$.
Proof: Let G be a smallest H-free with $\Delta \geq 9$ and $\chi > \max\{\omega, \Delta-1\}$. Then clearly as $G \neq C_{2n+1}$ or $K_{|V(G)|}$, $\chi = \Delta > \omega$. Let $u \in V(G)$. Then $G-u \neq K_{|V(G)|-1}$ (else $\chi = \omega$). If $\Delta(G-u) \geq 9$, then by minimality $\chi(G-u) \leq \max\{\omega(G-u), \Delta(G-u)-1\}$. Clearly if $\omega(G-u) \leq \Delta(G-u)-1$, then $\chi(G-u) = \Delta(G-u)-1 \leq \Delta-1$ and otherwise $\chi(G-u) = \omega(G-u) \leq \omega < \Delta$. In all the cases, $\chi(G-u) \leq \Delta-1$. Also if $\Delta(G-u) < 9$, then as $G-u \neq C_{2n+1}$ (else as G is H-free, G-u ~ $C_5$), by Brook's Theorem $\chi(G-u) \leq \Delta(G-u) < 9 \leq \Delta$. Thus always $\chi(G-u) \leq \Delta-1$ and in fact, $\chi(G-u) = \Delta-1$ and deg $u \geq \Delta-1$ $\forall u \in V(G)$.

Let $Q \subseteq V(G)$ be s.t. $<Q>$ is a maximum clique in G. Let $u \in Q$ be s.t. deg $u = \max_{v \in Q}$ deg v. Let $S = \{1,..., \Delta\}$ be a $\Delta$-coloring of V(G) s.t. u is colored $\Delta$ and $\{A_1, A_2,..., A_{\Delta-2}\} \subseteq N(u)$ where $A_i$ has color i for $1 \leq i \leq \Delta-2$. If $A_i$ is the only i-vertex in N(u), then $A_i$ is said to have a unique color.

**I.** Every $A_i$ with a unique color has at the most one repeat color and two vertices of that color (else $\overline{N(A_i)}$ has a color say r missing. Then color $A_i$ by r and u by i).
**II.** Now if deg $u = \Delta-1$, then all vertices in N(u) have unique colors (else some color r is missing in $\overline{N(u)}$. Color u by r) and if deg $u = \Delta$, then N(u) has $\Delta-2$ vertices with unique colors and two vertices with the same color. Thus N(u) has at least seven vertices with unique colors.
**III.** Every $A_i$ with a unique color has a j-vertex $\forall j \neq i$ (else color $A_i$ by j, u by i).

**Claim:** All vertices in N(u) with unique colors are adjacent to each other.
Let if possible say $A_r A_s \notin E(G)$ where $A_r$ and $A_s$ have unique colors. Then clearly $\exists A_r$-$A_s$ path = $\{A_r, B, C, A_s\}$ s.t. B, C have colors s, r resply. Now B is the only s-vertex in $V(G)-\overline{N(u)}$ (else let B' be another s-vertex. If $A_r B \notin E(G)$, then $<u, A_r, B, B'> = H$ and if $A_r B \in E(G)$, then $<A_r, B, B', A_s> =$

H). Similarly C is the unique r-vertex in V(G)-$\overline{N(u)}$. Hence $A_r$, $A_s$ are non-adjacent to at the most one more $A_k \in N(u)$, $k \notin \{r, s\}$ (else if say $A_rA_k$, $A_rA_l \notin E(G)$, then B has two s, k, l vertices and hence some color r missing in $\overline{N(B)}$. Color B by r, $A_r$ by s, u by r). By **II**, let $A_j$, $A_k$, $A_l$ be vertices in N(u) with unique colors s.t. $A_rA_i$, $A_sA_i \in E(G)$ for i = j, k, l. Also by **I**, w.l.g. let $A_j$ be the only j-vertex of $A_r$, $A_s$. Again by **I**, $A_j$ has either a unique r-vertex or s-vertex. W.l.g. let $A_r$ be the only r-vertex of $A_j$. Then color $A_j$ by r, $A_r$ and $A_s$ by j, u by s, a contradiction.

Hence the claim holds.

Thus deg u = $\Delta > \omega \geq \Delta-1 \Rightarrow \omega = \Delta-1$. Let Q = $\{u, A_1, A_2,..., A_{\omega-1}\}$ and $A_\omega$, $A_{\omega+1} \in N(u)$-Q have the same color $\omega$.

As G is H-free, V(G)-$\overline{N(u)}$ has no $\omega$-vertex. Hence $A_iA_\omega$ or $A_iA_{\omega+1} \in E(G)$ $\forall$ i, 1≤i≤$\omega$-1 (else color $A_i$ by $\omega$, u by i). W.l.g. let $A_\omega A_2 \notin E(G) \Rightarrow A_2A_{\omega+1} \in E(G)$ and let $A_{\omega+1}A_1 \notin E(G) \Rightarrow A_1A_\omega \in E(G)$. Now by **III**, $A_1$ ($A_2$) has an i-vertex $\forall$ i $\neq$ 1 (2). As $\Delta \geq 9$, w.l.g. let $A_3$, $A_4$, $A_5$ be adjacent to $A_\omega$. Again $A_\omega$ has at the most one repeat color (else $A_\omega$ has some color missing. Color $A_\omega$ by the missing color, $A_1$ by $\omega$, u by 1). Hence w.l.g. let $A_3$, $A_4$ be the only 3–vertex, 4-vertex of $A_\omega$. Also w.l.g. let $A_3$ be the only 3–vertex of $A_2$. Now if $A_3A_{\omega+1} \notin E(G)$, then color $A_3$ by $\omega$, $A_\omega$ by 3, $A_2$ by 3, u by 2 and if $A_3A_{\omega+1} \in E(G)$, then $A_1$ is the only 1–vertex of $A_3$ and hence color $A_3$ by 1, $A_1$ by $\omega$, $A_\omega$ by 3, $A_2$ by 3, u by 2, a contradiction in both the cases.

This proves Result 1.

**Main Result 2:** If $\Delta \geq 9$, and G is R-free, then $\chi \leq \max\{\omega, \Delta-1\}$.
Proof: Let G be a smallest R-free graph with $\Delta \geq 9$ and $\chi > \max\{\omega, \Delta-1\}$. Then clearly as G $\neq C_{2n+1}$ or $K_{|V(G)|}$, $\chi = \Delta > \omega$. Let u $\in$ V(G). Then G-u $\neq K_{|V(G)|-1}$ (else $\chi = \omega$). If $\Delta$(G-u) $\geq$ 9, then by minimality $\chi$(G-u) $\leq \max\{\omega$(G-u), $\Delta$(G-u)-1$\}$. Clearly if $\omega$(G-u) $\leq \Delta$(G-u)-1, then $\chi$(G-u) = $\Delta$(G-u)-1 $\leq \Delta$-1 and otherwise $\chi$(G-u) = $\omega$(G-u) $\leq \omega < \Delta$. In all the cases, $\chi$(G-u) $\leq \Delta$-1. Also if $\Delta$(G-u) < 9, then as G-u $\neq C_{2n+1}$ (else as G is R-free, G-u ~ $C_5$), by Brook's Theorem $\chi$(G-u) $\leq \Delta$(G-u) < 9 $\leq \Delta$. Thus always $\chi$(G-u) $\leq \Delta$-1 and in fact, $\chi$(G-u) = $\Delta$-1 and deg v $\geq \Delta$-1 $\forall$ v $\in$ V(G).

Let Q $\subseteq$ V(G) be s.t. <Q> is a maximum clique in G. Let u $\in$ Q be s.t. deg u = $\max_{v \in Q}$ deg v. Let S = $\{1,..., \Delta\}$ be a $\Delta$-coloring of V(G) s.t. u is colored $\Delta$ and $\{A_1, A_2,..., A_{\Delta-2}\} \subseteq$ N(u) where $A_i$ has color i for 1 $\leq$ i $\leq \Delta$-2. If $A_i$ is the only i-vertex in N(u), then $A_i$ is said to have a unique color.

**I.** Every vertex $A_i$ of N(u) with a unique color has at most two vertices of the same color (else $\overline{N(A_i)}$ has a color say r missing. Then color v by r and u by i).

**II.** Also as G is R-free, V(G)-$\overline{N(u)}$ has no repeat color (else if V, W $\in$ V(G)-$\overline{N(u)}$ have color say i, then <u, $A_i$, V, W> = R).
**III.** Every $A_i$ with a unique color has a j-vertex $\forall$ j $\neq$ i (else color $A_i$ by j, u by i).
**IV.** Now if deg u = $\Delta$-1, then all vertices in N(u) have unique colors (else some color r is missing in $\overline{N(u)}$. Color u by r) and if deg u = $\Delta$, then N(u) has $\Delta$-2 vertices with unique colors and two vertices with the same color. Thus N(u) has at least seven vertices with unique colors.

**Claim:** All vertices in N(u) with unique colors are adjacent to each other.
Let if possible say $A_1A_2 \notin E(G)$ where $A_1$, $A_2$ have unique colors. Then clearly $\exists A_1$-$A_2$ path = $\{A_1, B, C, A_2\}$ s.t. B, C have colors 2, 1 resply. Now by **II**, B is the only 2-vertex of $A_1$. Also $A_2$ is non-adjacent to at the most one more $A_k \in N(u)$, $k \neq 1$ (else if $A_2A_k$, $A_2A_l \notin E(G)$, then B has two 1, k, l vertices and hence has some color r missing in $\overline{N(B)}$. Color B by r, $A_1$ by 2, u by 1). By **IV**, w.l.g. let

$A_3,...,A_5$ be vertices in N(u) with unique colors s.t. $A_1A_i, A_2A_i \in E(G)$ for $3 \le i \le 5$. Also by **I**, w.l.g. let $A_3$ be the only 3 vertex of $A_1, A_2$. Again by **I**, $A_3$ has either a unique 1-vertex or 2-vertex. W.l.g. let $A_1$ be the only 1-vertex of $A_3$. Then color $A_3$ by 1, $A_1$ and $A_2$ by 3, u by 2, a contradiction.

Hence the claim holds.

Thus $\omega \ge \Delta-1 > \omega-1 \Rightarrow \omega = \Delta-1 \Rightarrow$ deg u = $\Delta$. Let Q = {u, $A_1, A_2,..., A_{\Delta-2}$} where $A_i$ has a unique color i and X, Y $\in$ N(u) be colored $\omega$. Further let $YA_1 \notin E(G)$. Then $XA_1 \in E(G)$ (else by **III**, $A_1$ has a $\omega$-vertex Z and $<A_1, Z, X, Y> = R$). Again w.l.g. let $XA_2 \notin E(G)$. Then $YA_2 \in E(G)$.

**Case 2.1:** V(G)-$\overline{N(u)}$ has no $\omega$-vertex.

By **III**, every $A_i$ is adjacent to X or Y. Also X (Y) is the only $\omega$-vertex of $A_1$ ($A_2$). Further as $\Delta \ge 9$, w.l.g. let $XA_i \in E(G)$ for $3 \le i \le 5$. Now X (Y) has a k-vertex $\forall k \ne \omega$ (else color X (Y) by the missing color, $A_1$ ($A_2$) by $\omega$, u by 1 (2)). Hence X and Y have at the most one repeat color. W.l.g. let $A_3$ be the unique 3-vertex of X and $A_2$. If $YA_3 \notin E(G)$, then color $A_3$ by $\omega$, X by 3, $A_2$ by 3, u by 2 and if $YA_3 \in E(G)$, then as $A_3$ has the unique 2-vertex $A_2$, color $A_3$ by 2, X by 3, $A_2$ by 3, $A_1$ by $\omega$, u by 1, a contradiction in both the cases.

**Case 2.2:** V(G)-$\overline{N(u)}$ has a $\omega$-vertex Z.

Then as G is R-free, every $A_j$ has two $\omega$-vertices and hence a unique j-vertex $A_j$, $j \ne i$.

**Case 2.2.1:** X has no r-vertex for some r, $2 \le r \le \omega-1$.

Then $XA_r \notin E(G) \Rightarrow YA_r \in E(G)$. Clearly $ZA_1, ZA_r \in E(G)$. Now Z has another 1-vertex (else color X by r, $A_1$ by $\omega$, Z by 1, u by 1). Also Y has an i-vertex $\forall i \ne \omega$ (else color X by r, Y by i, u by $\omega$). If $\exists$ s $\ne$ r, s.t. $YA_s \in E(G)$ and $A_s$ is the only s-vertex of Y, then color X by r, $A_s$ by 1, $A_1$ by s, Y by s, u by $\omega$, a contradiction. Hence whenever $YA_s \in E(G)$, s $\ne$ r, Y has another s-vertex. Hence Y is adjacent to at the most one $A_s$, s $\ne$ r. As $\Delta \ge 9$, w.l.g. let $XA_j, ZA_j \in E(G)$ for say i = 3 $\le$ i $\le$ 6. As Z has two 1-vertices, Z has at the most one more repeat color (else $\overline{N(Z)}$ has some color t missing, color X by r, Z by t, $A_1$ by $\omega$, u by 1) and hence w.l.g. let Z have the unique 3-vertex $A_3$. Color X by r, Z by 3, $A_3$ by $\omega$, u by 3, a contradiction.

**Case 2.2.2**: X, Y have an i-vertex $\forall i \ne \omega$.

Clearly $ZA_1, ZA_2 \in E(G)$. As $\Delta \ge 9$, w.l.g. $XA_j \in E(G)$ for i = 3, 4, 5. Also as X has at the most one repeat color, w.l.g. let $A_j$ be the only i-vertex of X for i = 3, 4. Now Z has another 1-vertex (else color $A_3$ by 2, $A_2$ by 3, X by 3, $A_1$ by $\omega$, Z by 1, u by 1) and Z has no color missing (else color Z by the missing color t, $A_3$ by 2, $A_2$ by 3, X by 3, $A_1$ by $\omega$, u by 1). Hence Z has at the most one more repeat color other than 1 $\Rightarrow ZA_3 \notin E(G)$ or $ZA_4 \notin E(G)$ (else Z has a unique say 4-vertex $A_4$. Color X by 3, $A_3$ by 2, $A_2$ by 3, $A_4$ by $\omega$, Z by 4, u by 4). W.l.g. let $ZA_3 \notin E(G) \Rightarrow YA_3 \in E(G)$. If $A_2$ is the unique 2-vetex of Y, then color $A_1$ by 2, $A_2$ by 3, $A_3$ by 1, X by 3, Y by 2, u by $\omega$ and if $A_3$ is the unique 3-vetex of Y, then color $A_3$ by 1, $A_1$ by 3, Y by 3, $A_4$ by 2, $A_2$ by 4, X by 4, u by $\omega$, a contradiction in both the cases.

This proves Result 2.